\documentclass{article}
\usepackage[none]{hyphenat}
\usepackage{xcolor}
\usepackage{amsmath,amsfonts,graphicx,amsthm,geometry,setspace,tikz-cd,subcaption,url}
\usepackage{etoolbox, authblk}
 \newtheorem{thm}{Theorem}[section]
 \newtheorem*{theorem*}{Main Theorem}
 \newtheorem*{thm*}{Theorem}
 \newtheorem{cor}[thm]{Corollary}
 
 \newtheorem{prop}[thm]{Proposition}
 \theoremstyle{definition}
 \newtheorem{defn}[thm]{Definition}
 \theoremstyle{remark}
 \newtheorem{rem}[thm]{Remark}
 
 \numberwithin{equation}{section}

\title{Virtual Knotoids in Thickened Surfaces}
\author[1]{Neslihan G\"ug\"umc\"u}
\author[1]{Hamdi Kayaslan}
\affil[1]{Department of Mathematics, Izmir Institute of Technology, Urla Izmir 35430, Turkey}

\begin{document}

\maketitle

\begin{abstract}
In this paper, we give an interpretation of virtual knotoids as rail arcs that are embedded in thickened surfaces. We prove that virtual knotoids admit a unique irreducible rail arc representation in thickened surfaces. We also prove a conjecture of Kauffman and the first author stating that virtual knotoid theory is a proper generalization of classical knotoid
theory. 
\end{abstract}

\section{Introduction}
Knotoids are open-ended knot diagrams that lie on a surface. The theory of knotoids was introduced by Vladimir Turaev \cite{turaev2012knotoids} in 2012, as a theory of generalized open-ended knot diagrams or $(1, 1)$-tangles. The theory of knotoids and its applications in topological analysis of proteins have been further studied by many researchers \cite{adams2019knots,goundaroulis2017studies,goundaroulis2017topological,gugumcu2022invariants,gugumcu2017knotoids,moltmaker2022framed}, including the first author and Louis H. Kauffman \cite{gugumcu2017new,gugumcu2021parity,gugumcu2021quantum}.

As initially remarked in \cite{turaev2012knotoids}, it is natural to examine knotoids in the context of virtual knot theory. Virtual knot theory was introduced by Louis H. Kauffman \cite{kauffman2012introduction} as a generalization of classical knot theory, which studies knots in three-dimensional Euclidean space. 

\begin{figure}[h!]
    \centering
    \includegraphics[width=0.45\linewidth]{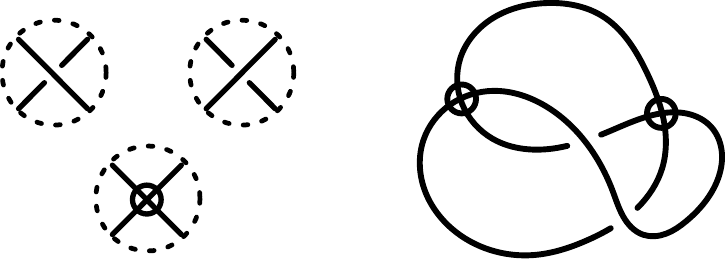}
    \caption{Classical crossings and a virtual crossing on the left-hand side. A virtual knot diagram with two classical and two virtual crossings on the right-hand side.}
    \label{vk}
\end{figure}

A virtual knot (or link) diagram may contain two types of crossings: classical and virtual. See Figure \ref{vk}. Virtuality of a crossing refers to non-planarity of the graph underlying the diagram. More precisely, a virtual knot diagram can be assigned to a knot in some thickened surface, studied up to ambient isotopy in the thickened surface and addition or deletion of hollow handles to or from the thickened surface. 

 A representation of a virtual knot in a thickened surface is called \textit{irreducible} if there is no hollow handle to cut off and reduce the genus of the thickened surface. Greg Kuperberg showed in \cite{kuperberg2003virtual} that every virtual link admits a unique irreducible representation in thickened surfaces up to isotopy of thickened surfaces.

In \cite{gugumcu2017new}, virtual knotoids were introduced diagrammatically as open-ended virtual knot diagrams in $S^2$ and was also shown in \cite{gugumcu2017new} that the theory of knotoids in surfaces up to ambient isotopy of surfaces and addition or deletion of hollow handles is equivalent to the theory of virtual knotoids. 

In this paper we give a three-dimensional interpretation of virtual knotoids with \textit{rail arcs}. A rail arc is an open-ended smooth curve in a thickened surface. The endpoints of a rail arc are attached on two vertical parallel line segments that we call \textit{rails}. We define \textit{arc isotopy} of rail arcs as an ambient isotopy of the thickened surface that turns a rail arc into another in the complement of two rails. Rail arcs are considered up to arc isotopy plus addition or deletion of hollow handles. We first prove that a rail arc in a thickened surface can be represented as a knotoid diagram in the underlying surface subject to Reidemeister moves that take place away from the endpoints of the knotoid diagram. It follows from this correspondence and the observation of Kauffman and the first author on surface knotoids and virtual knotoids that rail arcs can be uniquely represented by virtual knotoids. 

Due to the nature of the equivalence relation assumed for rail arcs, a rail arc can be represented in thickened surfaces of different genera, and it is natural to ask if the minimum genus representation of a rail arc is unique. We give an affirmative answer to this question by proving the following theorem.

\begin{theorem*}\label
{thm:main}
	Every rail arc has a unique irreducible representation up to arc isotopy.
\end{theorem*} 

We prove this theorem by extending Kuperberg's arguments in \cite{kuperberg2003virtual} to the case of rail arcs. We give an illustrative exposition of the arguments from 3-manifold theory and study how the arguments extend in the presence of endpoints and rails to which they are attached.

Furthermore, Kauffman and the first author conjectured in \cite{gugumcu2017new} that classical knotoid theory (that is, the theory of knotoids in the plane with only classical crossings that are considered up to classical Reidemeister moves) embeds properly into the theory of virtual knotoids. We give a proof of this conjecture by the help of the Main Theorem.

Let us now give an outline of the paper. In Section \ref{review} we review the basics of classical and virtual knotoids. In particular, in Section \ref{knotoids}, we introduce knotoid diagrams in a surface $\Sigma$, their equivalence in $\Sigma$, and the semigroup structure of the equivalence classes of knotoid diagrams in $S^2$. We discuss the relation between classical knots and knotoids. We mention geometric interpretations of knotoids in $S^2$ \cite{turaev2012knotoids} and knotoids in $\mathbb{R}^2$ \cite{gugumcu2017new}. In Section \ref{vknotoids}, we review virtual knotoids and their relations with knotoids in surfaces.  In Section \ref{virtualarcs}, we establish a three-dimensional interpretation of virtual knotoids in $S^2$. To do this, we first introduce rail arcs in thickened surfaces in Section \ref{varcs} and then show that rail arcs are in one-to-one correspondence with knotoids in surfaces. In Section \ref{main}, we study virtual knotoids as rail arcs in thickened surfaces and prove the Main Theorem. 

\section{A review on knotoids}
\label{review}

We begin with a review on fundamental notions of the theory of knotoids.

\subsection{Knotoids}
\label{knotoids}
\begin{figure}[h!]
	\centering
	\includegraphics[width=0.65\linewidth]{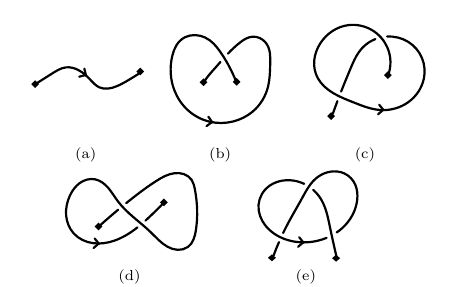}
	\caption{Examples of knotoid diagrams.}
	\label{kdiagrams}
\end{figure}

A knotoid diagram is an open-ended knot diagram in an oriented surface. Examples of knotoid diagrams are shown in Figure \ref{kdiagrams}. A more formal definition of a knotoid diagram is as follows.
\begin{defn}
	A \textit{ knotoid diagram} is a generic immersion of the unit interval $[0,1]$ in the interior of a surface $\Sigma$ with finitely many transversal double self-intersections. These self-intersections are endowed with over or under passage information, and together with this information they are called \textit{crossings} of the knotoid diagram. The images of $0$ and $1$ are two distinct points that are considered to be the \textit{endpoints} of the knotoid diagram and are called the \textit{tail} and the \textit{head} of the diagram, respectively. A knotoid diagram is considered to be oriented from its tail to its head. A \textit{trivial knotoid diagram} is an embedding of the unit interval in $\Sigma$, that is, it is a knotoid diagram with no crossings as shown in Figure \ref{kdiagrams}a.
\end{defn}

\begin{figure}[h!]
	\centering    
	\begin{subfigure}{0.5\textwidth}
		\centering
		\includegraphics[width=38mm]{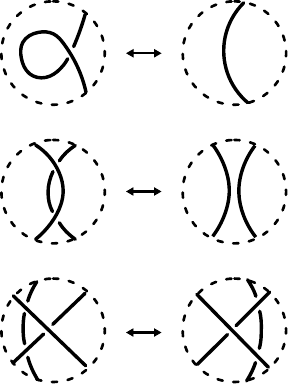}
		\caption{Reidemeister moves.}
		\label{reidemoves}
	\end{subfigure}\hfill
	\begin{subfigure}{0.5\textwidth}
		\centering
		\includegraphics[width=45mm]{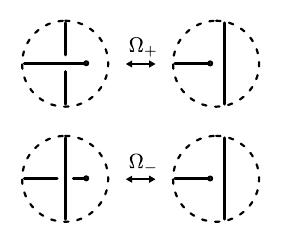}
		\vspace{1 cm}
		\caption{Forbidden $\Omega_+$ and $\Omega_-$ moves.}
		\label{forbidden}
	\end{subfigure}
	\caption{Moves for knotoids.}
	\label{knotoidmoves}
\end{figure}

Knotoid diagrams in a surface $\Sigma$ are considered up to the equivalence relation induced by three \textit{Reidemeister moves}, depicted in Figure \ref{reidemoves}, and isotopy of $\Sigma$. All of these moves change a knotoid diagram locally, as shown in Figure \ref{reidemoves}, and are assumed to be performed away from the endpoints of the knotoid diagram. This means that it is not allowed to pull or push a strand with an endpoint over or under another strand as depicted in Figure \ref{forbidden}. We denote these forbidden moves by $\Omega_+$ and $\Omega_-$, respectively. Notice that any knotoid diagram would be equivalent to the trivial knotoid diagram if $\Omega_+$ and $\Omega_-$ moves were allowed.

\begin{defn}
	A \textit{knotoid} in a surface $\Sigma$ is an equivalence class of knotoid diagrams, with respect to the equivalence relation induced by Reidemeister moves and isotopy of $\Sigma$. The set of knotoids in $\Sigma$ is denoted by $\mathcal{K}(\Sigma)$. In particular, knotoids in $S^2$ are called \textit{spherical knotoids}, and knotoids in $\mathbb{R}^2$ are called \textit{planar knotoids}.
\end{defn}

\begin{figure}[h!]
	\centering
	\includegraphics[width=0.5\linewidth]{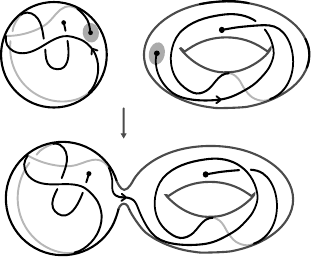}
	\caption{Product of a knotoid diagram in sphere and a knotoid diagram in torus.}
	\label{consum}
\end{figure}

One can define a product operation of two knotoids in surfaces as follows.
\begin{defn}
	Let $\kappa_i$ be a knotoid in $\Sigma^i$, and $K_i$ be a diagram in the equivalence class $\kappa_i$ having endpoints  $t_i$ and  $h_i$ for $i=1,2$. Let $U$ and $V$ be regular neighborhoods of $h_1$ and $t_2$, respectively. The \textit{product} of knotoids $\kappa_1$ and $\kappa_2$, denoted by $\kappa_1\kappa_2$, is the equivalence class of the knotoid diagram resulting from the gluing of $\Sigma^1-\text{Int}(U)$ and $\Sigma^2-\text{Int}(V)$ along a homeomorphism $\phi:\partial(U)\rightarrow \partial(V)$ such that the intersection point $K_1\cap \partial(U)$ gets mapped to the intersection point $K_2 \cap \partial(V)$. $K_1 K_2$ is a knotoid diagram in $\Sigma=\Sigma^1\#\Sigma^2$ having tail $t_1$ and head $h_2$, see Figure \ref{consum}. 
\end{defn}

Note that the set of knotoids in $S^2$ endowed with the product operation forms a semigroup with the identity element that is the trivial knotoid.

\begin{rem}
	The notion of knotoids can be extended by allowing the domain of immersion to include a finite number of circles in addition to the unit interval. The resulting diagrams are called \textit{multi-knotoid diagrams.} \textit{Multi-knotoids} are then equivalence classes of multi-knotoid diagrams considered up to Reidemeister moves and isotopy of $\Sigma$.
\end{rem}

Every planar knotoid diagram can be mapped to a spherical knotoid diagram by adding the point at infinity to $\mathbb{R}^2$. However, there exist non-trivial planar knotoid diagrams which map to spherical knotoid diagrams that are equivalent to trivial knotoid diagrams in $S^2$. An example is shown in Figure \ref{kdiagrams}b. The reason behind this is as follows. Think of $S^2$ as $\mathbb{R}^2\cup \infty$. Then it is clear that the spherical isotopy not only includes planar isotopy but also lets an arc of a diagram pass through the point at infinity, which transforms the non-trivial planar knotoid diagram in Figure \ref{kdiagrams}b into a trivial one in $S^2$. This implies that planar knotoids are not in a one-to-one correspondence with spherical knotoids.

\begin{figure}[h!]
	\centering
	\includegraphics[width=0.4\linewidth]{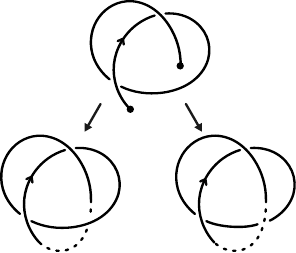}
	\caption{The underpass (on the left) and overpass (on the right) closures of a knotoid diagram.} 
	\label{uoclosure}
\end{figure}

One can tie up the endpoints of a knotoid diagram in $S^2$ or $\mathbb{R}^2$ to obtain a (classical) knot diagram with \textit{underpass} or \textit{overpass closures}. The underpass (overpass) closure of a knotoid diagram is obtained by connecting the endpoints of a knotoid diagram by an arc which goes under (over) every other arc it meets along the way. The resulting diagram is clearly a knot diagram which represents a knot in $S^3$ (a classical knot). The underpass and overpass closures may represent inequivalent knots. We exemplify this in Figure \ref{uoclosure} where the underpass closure of a knotoid diagram represents the trefoil knot, while the overpass closure represents the unknot. To have a well-defined closure on knotoid diagrams, we fix the closure type we want to use to obtain knots from knotoids. We have the following theorem by fixing the closure as the underpass closure.

\begin{thm}
	\cite{turaev2012knotoids}\label{thm:knotoidrepresetantion} The underpass closures of two knotoid diagrams in $S^2$ represent the same knot if and only if they are related to each other by finitely many Reidemeister moves and forbidden $\Omega_-$ moves.
\end{thm}

By Theorem \ref{thm:knotoidrepresetantion} we can consider the collection of spherical knotoid diagrams representing a classical knot via the underpass closure as \textit{knotoid representations} of the knot. In \cite{turaev2012knotoids}, Turaev suggested to utilize knotoids in computing knot invariants, as they provide simpler representations of classical knots possibly with fewer crossings than knot diagrams.

A knotoid in $S^2$ that can be represented by a knotoid diagram with endpoints lying in the same local region of the sphere determined by the diagram  topologically coincide with the classical knot obtained by closing its endpoints. Such knotoid diagrams are called \textit{knot-type knotoid diagrams}. There exists another type of spherical knotoids that we call \textit{proper knotoids}. A proper knotoid is a knotoid without any representative knotoid diagrams in its equivalence class that has its endpoints in the same region. In this manner, the spherical knotoid theory is a proper extension of the classical knot theory. 

In \cite{turaev2012knotoids}, Turaev showed that spherical knotoids can be interpreted as \textit{theta-curves} in $S^3 = \mathbb{R}^3 \cup \infty$.

\begin{defn}
	A \textit{theta-curve} is a directed graph embedded in $S^3$ consisting of two vertices, labeled $v_{0}, v_{1}$ and three edges connecting these two vertices labeled $e_{+}, e_{0}, e_{-}$, considered up to label-preserving isotopies of $S^3$. The edges of a theta-curve is directed from $v_0$ to $v_1$.  A theta-curve is called \textit{simple} if the edges $e_+$ and $e_-$ bound a unique disk in $S^3$.
\end{defn}

\begin{figure}[h!]
	\centering    
	\begin{subfigure}{0.75\textwidth}
		\centering
		\includegraphics[width=76mm]{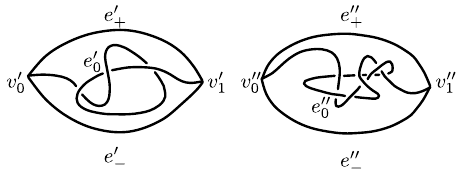}
		\vspace{0.1 cm}
		\caption{Two theta-curve examples, $\theta'$ and $\theta''$.}
		\label{thtex}
	\end{subfigure}\hfill
	\begin{subfigure}
		{0.72\textwidth}
		\vspace{0.3cm}
		\centering
		\includegraphics[width=67mm]{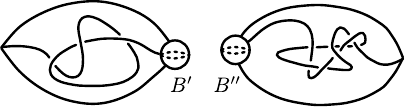}
        \vspace{0.2 cm}
		\caption{Two regular neighbourhoods $B'$ and $B''$ around the vertices $v_1'$ and $v_0''$.}
		\label{balls}
	\end{subfigure}
	\begin{subfigure}
		{0.77\textwidth}
		\vspace{0.1cm}
		\centering
		\includegraphics[width=78mm]{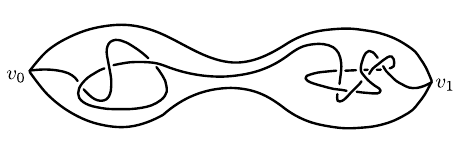}
		\vspace{-0.2cm}\caption{$\theta=\theta'\theta''.$}
		\label{thtmult}
	\end{subfigure}
	\caption{The vertex multiplication of theta-curves. }
	\label{tht}
\end{figure}

\begin{defn}
	Let $\theta'$ and $\theta''$ be two theta-curves. Let $B'$ and $B''$ be regular neighborhoods of vertices $v_1'$ of $\theta'$ and $v_0''$ of $\theta''$, respectively. The \textit{vertex multiplication} of $\theta'$ and $\theta''$, denoted by $\theta'\theta''$, is a theta-curve resulting from the gluing of $S^3-\text{Int}(B')$ and $S^3-\text{Int}(B'')$ along an orientation-reversing homeomorphism $\varphi:\partial(B')\rightarrow \partial(B'')$ such that the intersection point $e_i'\cap\partial(B')$ gets mapped to the intersection point $e_i''\cap\partial(B'')$ for $i=+,0,-$. See Figure \ref{tht} for an illustration of vertex multiplication. The vertex multiplication on theta-curves induces a semigroup structure on theta-curves.
\end{defn}

One can construct a theta-curve $\theta$ from a knotoid diagram $K$ in $\mathbb{R}^2$. Consider $\mathbb{R}^2$ as $\mathbb{R}^2\times \{0\}\subset\mathbb{R}^3$, and the tail and the head of $K$ as the vertices $v_0$ and $v_1$ of $\theta$, respectively. The edge $e_0$ of $\theta$ is formed by pushing the overpassing strands of $K$ in $\mathbb{R}^2\times \mathbb{R}^+$. The edges $e_{+}$ and $e_{-}$ of $\theta$ are formed by pushing an embedded arc in $\mathbb{R}^2\times \{0\}$ connecting the vertices $v_0$ and $v_1$ to $\mathbb{R}^2\times\mathbb{R}^+$ and $\mathbb{R}^2\times\mathbb{R}^-$, respectively. It is clear that $e_+ \cup e_0 \cup e_-$ is a simple theta-curve in $S^3$.

This construction induces a bijective map between spherical knotoids and simple theta-curves in $S^3$. In fact, we have the following theorem.
\begin{thm}\cite{turaev2012knotoids}
	There is an isomorphism between the semigroup of knotoids in $S^2$ and the semigroup of simple theta-curves in $S^3$.
\end{thm}

\begin{figure}[h!]
	\centering
	\includegraphics[width=110mm]{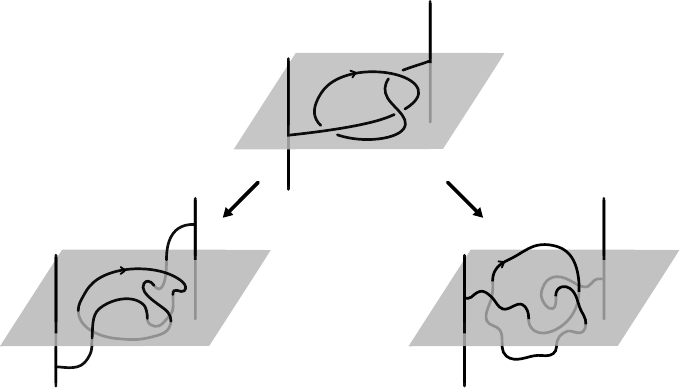}
	\caption{Two different embeddings of a knotoid diagram in $\mathbb{R}^3$.}
	\label{interp}
\end{figure}

In \cite{gugumcu2017new}, the first author and Kauffman gave a three-dimensional interpretation for planar knotoids as follows. Let $K$ be a knotoid in $\mathbb{R}^2$. Identify $\mathbb{R}^3$ with $\mathbb{R}^2\times \mathbb{R}$ so that $K$ can be considered lying in $\mathbb{R}^2\times \{0\}$, and can be embedded in $\mathbb{R}^3$ by pushing the overpassing strands to $\mathbb{R}^2\times \mathbb{R^+}$ and underpassing strands to $\mathbb{R}^2\times \mathbb{R^-}$, while the endpoints of $K$ are kept attached on two distinct parallel lines that are orthogonal to the plane where $K$ lies. In this way, an open, oriented smooth curve embedded in $\mathbb{R}^3$ is obtained. See Figure \ref{interp} for an illustration of such embeddings.

Conversely, any open, oriented curve in $\mathbb{R}^3$ that is generic to the \textit{xy}-plane can be projected onto the plane to obtain a knotoid diagram. The one-to-one correspondence between knotoids in $\mathbb{R}^2$, and open, oriented smooth curves in $\mathbb{R}^3$ is obtained by defining \textit{line isotopy} for space curves.

\begin{defn}
	Let $K_1$ and $K_2$ be two smooth, open and oriented curves embedded in $\mathbb{R}^3$ with the endpoints corresponding to points on the lines $t_1\times \mathbb{R}$, $h_1\times\mathbb{R}$ and $t_2\times \mathbb{R}$, $h_2\times\mathbb{R}$. $K_1$ and $K_2$ are said to be \textit{line isotopic} if there is a smooth ambient isotopy of $\mathbb{R}^3$, taking one curve to the other curve in the complement of the lines, taking the tail and the head of $K_1$ to the tail and the head of $K_2$, and taking lines to lines; $t_1 \times \mathbb{R}$ to $t_2 \times \mathbb{R}$ and $h_1 \times \mathbb{R}$ to $h_2 \times \mathbb{R}$.
	\label{li}
\end{defn}

\begin{thm} \cite{gugumcu2017new} Two open oriented curves embedded in $\mathbb{R}^3$, which are both generic to  \textit{xy}-plane, are line isotopic (with respect to the lines determined by the endpoints of the curves) if and only if the projections of the curves to $\mathbb{R}^2\times \{0\}$ are equivalent knotoid diagrams.
\end{thm}

The approach used in the geometric interpretation of planar knotoids given above has found key applications in studying protein folding. See \cite{goundaroulis2017studies,goundaroulis2017topological,gugumcu2022invariants}.

\subsection{Virtual knots and virtual knotoids}
\label{vknotoids}
In this section, we review the basic notions of virtual knots and virtual knotoids. The theory of virtual knots was introduced by Louis H. Kauffman in his paper \cite{kauffman1999virtual}. 
Although virtual knots are realized as knots in thickened surfaces modulo an extended equivalence relation, they can equivalently be studied through their diagrams in the plane.

\begin{figure}[h!]
    \centering
    \includegraphics[width=0.09\linewidth]{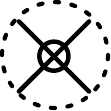}
    \caption{A Virtual crossing.}
    \label{vcross}
\end{figure}

We begin with a diagrammatic description of virtual knots and later present a topological description.
\begin{defn}
\label{fvk}
	A \textit{virtual knot diagram} is an immersion of a circle in the plane or 2-sphere $S^2$ whose each self-intersection is declared as either a classical crossing or a \textit{virtual crossing}. A virtual crossing is indicated by a circle around a transversal intersection that provides no under or over information unlike a classical crossing. See Figure \ref{vcross}. A virtual knot diagram without any virtual crossings is called a \textit{classical knot diagram}.

\end{defn}

\begin{figure}[h!]
	\centering    
	\begin{subfigure}{0.5\textwidth}
		\centering
		\includegraphics[width=38mm]{classicalreidemoves2}
		\caption{Classical Reidemeister moves.}
		\label{vcrm}
		\vspace{0.3 cm}
	\end{subfigure}\hfill
	\begin{subfigure}{0.5\textwidth}
		\centering
		\includegraphics[width=38mm]{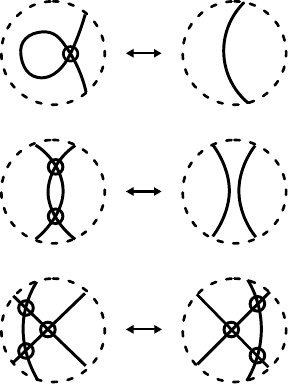}
		\caption{Virtual Reidemeister moves.}
		\vspace{0.3 cm}
		\label{vvrm}
	\end{subfigure}
	\begin{subfigure}{0.5\textwidth}
		\centering
		\includegraphics[width=38mm]{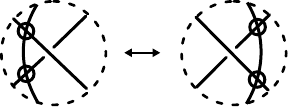}
		\vspace{0.3 cm}
		\caption{Mixed virtual move.}
		\label{vmm}
	\end{subfigure}\hfill
	\caption{Moves for virtual knot diagrams in the plane.}
	\label{virtualknotoidmoves}
\end{figure}

Virtual knot diagrams are considered up to an equivalence relation called the \textit{virtual equivalence}. The virtual equivalence is induced by the classical Reidemeister moves shown in Figure \ref{vcrm}, plus the detour moves. A \textit{detour move} allows an arc containing a sequence of only virtual crossings to
be replaced by another arc that contains a (possibly empty) sequence of only virtual
crossings.
Some of the special detour moves that resemble classical Reidemeister moves are depicted in Figures \ref{vvrm} and \ref{vmm}. Notice that we call the move in Figure \ref{vmm} a mixed virtual move since it involves one classical crossing and two virtual crossings. In this move, the detour move is applied to the arc that contains virtual crossings and as a result this arc is pulled across the classical crossing.

\begin{defn}
	A \textit{virtual knot} is an equivalence class of virtual knot diagrams in the plane up to the virtual equivalence. 
\end{defn}

\begin{rem}
Virtual knots can be extended to have more components. We call a union of virtual knots a \textit{virtual link.} Virtual links are subject to the same equivalence relation defined on virtual knots.
\end{rem}

\begin{figure}[h!]
    \centering
    \includegraphics[width=0.4\linewidth]{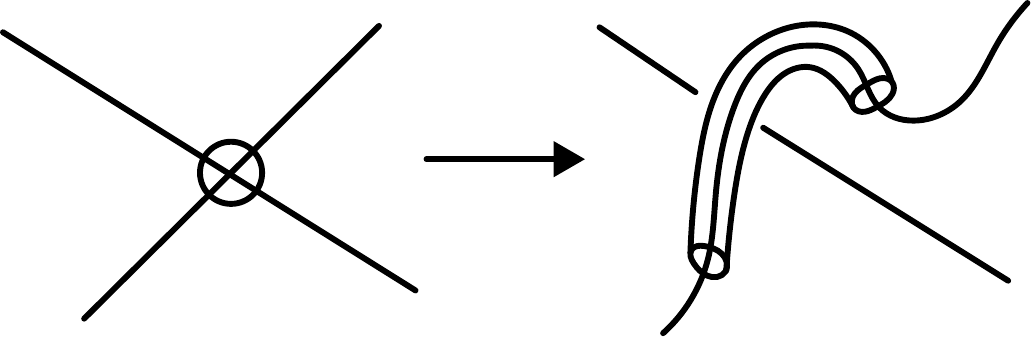}
    \caption{One of the arcs in a virtual crossing is held by a 1-handle.}
    \label{handle}
\end{figure}

Let $K$ be a virtual knot diagram and $I$ be the unit interval $[0,1]$. One can attach a 1-handle $S^1\times I$ to $S^2$ for each virtual crossing of $K$ so that one of the arcs at the virtual crossing passes through the handle, as depicted in Figure \ref{handle}. In this way $K$ is associated a knot diagram without virtual crossings in a closed, connected and orientable surface. The resulting knot diagram in the surface is called a \textit{surface representation} of $K$. 

In the sequel, a surface is always assumed to be closed, connected and orientable and it is denoted by $\Sigma_g$ where $g$ is the genus of the surface.

\begin{defn}
	Let $K$ be a knot or link diagram in a surface $\Sigma_g$. Let $B_1, B_2$ be two disjoint open disks in $\Sigma_g$ such that $B_1$ and $B_2$ do not intersect $K$. The \textit{stabilization} operation consists of removing $B_1$ and $B_2$ from $\Sigma_g$, and attaching a 1-handle $S^1\times I$ by gluing the circles $S^1\times \{0\}$ and $S^1\times \{1\}$ to the boundaries of removed disks in $\Sigma_g$. 
    Let $C$ be a simple closed curve in $\Sigma_g$ such that $C$ is disjoint from $K$ and does not bound a disk in $\Sigma_g-K$. The \textit{destabilization} operation consists of cutting $\Sigma_g$ along $C$, and capping the two resulting components with disks $D^2$.
\end{defn}

\begin{defn}
    The equivalence relation on knot diagrams in surfaces induced by Reidemeister moves, homeomorphisms of surfaces and handle stabilization and destabilization operations of surfaces is called \textit{stable equivalence}. Two knot diagrams in surfaces that are equivalent up to stable equivalence are called \textit{stably equivalent}.
\end{defn}

\begin{thm}\label{thmx}\cite{kamada2000abstract}
  Two virtual knot diagrams are virtually equivalent if and only if their surface representations are stably equivalent.
\end{thm}

A knot diagram $K$ in a surface $\Sigma_g$ can be lifted to a knot that is embedded in the thickened surface $\Sigma_g \times I$, where $I$ denotes the unit interval $[0,1]$ in the following way. We identify $\Sigma_g$ with $\Sigma_g\times \{a\}$, where $a\in (0,1)$. We push the overpassing arcs of $K$ into $\Sigma_g\times (a,1)$ and underpassing arcs into $\Sigma\times (0,a)$. The stabilization and destabilization operations defined on surfaces are extended to thickened surfaces to study such embeddings as follows. 

\begin{figure}[h!]
	\centering
	\includegraphics[width=0.55\linewidth]{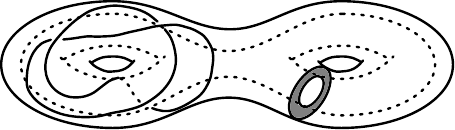}
	\caption{A knot in a thickened genus two surface and a destabilization annulus.}
	\label{desanknot}
\end{figure}

\begin{defn} Let $V$ be an embedding of a circle in a thickened surface $\Sigma_g \times I$. Let $B_1, B_2$ be two disjoint open disks in $\Sigma_g$ such that $B_1\times I$ and $B_2\times I$ do not intersect $V$. The \textit{stabilization} operation consists of removing the two open annuli $B_1\times I$ and $B_2\times I$ from $\Sigma_g\times I$, and attaching a thickened handle $S^1\times I \times I$ by gluing the annuli $S^1\times \{0\}\times I$ and $S^1\times \{1\}\times I$ to the boundaries of removed annuli in $\Sigma_g\times I$. Let $C$ be a simple closed curve in $\Sigma_g$ such that $C\times I$ is disjoint from $V$ and does not bound a ball in $\Sigma_g\times I-V$. The \textit{destabilization} operation consists of cutting $\Sigma_g\times I$ along the vertical annulus $C\times I$, and capping the two resulting annuli with thickened disks $D^2\times I$. The annulus $C\times I$ is called a $\textit{destabilization annulus}$. Figure \ref{desanknot} illustrates a destabilization annulus.
	\label{destab}
\end{defn}

\begin{defn}
    The equivalence relation on knots in thickened surfaces induced by isotopies of thickened surfaces and handle stabilization and destabilization operations of thickened surfaces is called \textit{stable equivalence}. Two equivalent knots in thickened surfaces up to stable equivalence are called \textit{stably equivalent}.
\end{defn}

By extending the stable equivalence to knots in thickened surfaces, we have the following results for virtual knots.

\begin{prop}
Two knots in thickened surfaces $\Sigma_{g_1} \times I$ and $\Sigma_{g_2} \times I$ are stably equivalent if and only if their surface representations in $\Sigma_{g_1}$ and $\Sigma_{g_2}$ are stably equivalent. 
    
\end{prop}

\begin{thm}\label{thmy}\cite{kauffman1999virtual,kuperberg2003virtual}
Two virtual knot diagrams are virtually equivalent if and only if the corresponding knots in thickened surfaces are stably equivalent.
\end{thm}

\begin{figure}[h!]
	\centering
	\includegraphics[width=0.6\linewidth]{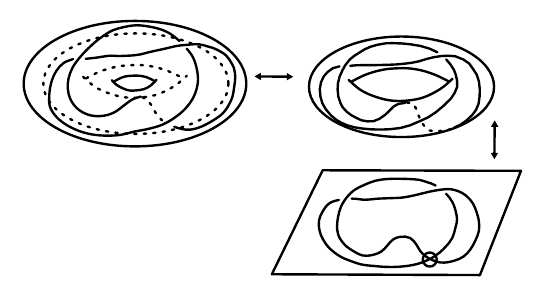}
    \vspace{-0.5 cm}
	\caption{A virtual knot with virtual genus one and its representations in torus and thickened torus.}
	\label{vknots}
\end{figure}

It is clear that stabilization and destabilization operations on surfaces change the genus of a surface so that a virtual knot can be represented by knots lying in surfaces of various genera. A natural problem arising from this fact is finding the minimal genus among the surfaces that a virtual knot can be represented in. The minimum genus of a surface in which a virtual knot $K$ can be represented, is called the \textit{virtual genus} of the virtual knot $K$. A virtual knot with virtual genus one is shown in Figure \ref{vknots}. The virtual genus is an invariant of virtual knots and was studied in \cite{boden2021minimal, chrisman2025uq, dye2005minimal, kamada2000abstract}.

\begin{defn}
    A representation of a virtual knot in $\Sigma_g \times I$ is called \textit{irreducible} if $g$ is equal to the virtual genus of the virtual knot.
\end{defn}     
    
Another natural question concerning virtual genus is whether a representation of a virtual knot in a thickened surface of minimal genus is unique up to isotopy of that thickened surface. In \cite{kuperberg2003virtual}, Kuperberg gave an affirmative answer to this question by proving Theorem \ref{kup} for virtual links. In Section \ref{main}, we adapt this question  for virtual knotoids and prove our main result. (See Theorem \ref{thm1}).

\begin{thm}
	\label{thm:Kuperberg}
	\cite{kuperberg2003virtual} Every virtual link has a unique irreducible representation in thickened surfaces up to isotopy.
    \label{kup}
\end{thm}

In \cite{turaev2012knotoids}, it was mentioned that knotoids can be studied in virtual setting. This can be done by obtaining a virtual knot diagram from a knotoid diagram in $\mathbb{R}^2$ or $S^2$ using the \textit{virtual closure} or considering a knotoid diagram with virtual crossings, as we explain in the following.

\begin{figure}[h!]
	\centering
	\includegraphics[width=0.4\linewidth]{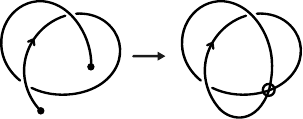}
	\caption{A virtual closure of a knotoid diagram.}
	\label{vclo}
\end{figure}

\begin{defn}
	A \textit{virtual closure} of a knotoid diagram $K$ connects the tail and head of $K$ by an arc so that each of its intersections with $K$ is declared as a virtual crossing (see Figure \ref{vclo}). 
\end{defn}

In \cite{gugumcu2017new}, it was shown that virtual closure induces a well-defined map from the set of spherical knotoids to the set of virtual knots. Being a well-defined map, any virtual knot invariant can be considered as an invariant of spherical knotoids through the virtual closure.  See \cite{gugumcu2017new} for more details. 

In the same paper, the first author and Kauffman also introduced virtual knotoids. Let us now review basic notions of virtual knotoids.

\begin{figure}[h!]
	\centering
	\includegraphics[width=0.6 \linewidth]{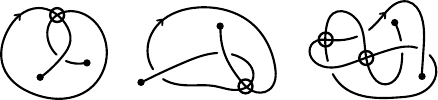}
	\caption{Virtual knotoid diagrams.}
	\label{virtk1}
\end{figure}

\begin{defn}
	A \textit{virtual knotoid diagram} is a knotoid diagram in $S^2$ that may also contain a finite number of virtual crossings. Some examples of virtual knotoid diagrams are given in Figure \ref{virtk1}. Specifically, a virtual knotoid diagram that contains only classical crossings is called a \textit{classical knotoid diagram}.
\end{defn}

\begin{figure}[h!]
	\centering        \includegraphics[width=45mm]{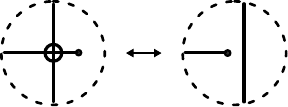}
	\vspace{0.3 cm}
	\caption{Virtual endpoint move.}
	\label{vem}
\end{figure}

Two virtual knotoid diagrams are said to be \textit{virtually equivalent} if they are related to each other by a sequence of classical Reidemeister moves and the detour moves as virtual knot diagrams. Detour moves induce the \textit{virtual endpoint move} depicted in Figure \ref{vem}, which is to pull the arc adjacent to an endpoint  through a transversal arc by deleting or adding a virtual crossing. Note that it is forbidden to pull an endpoint of a virtual knotoid diagram over or under a strand since such moves transform any virtual knotoid diagram into the trivial knotoid diagram without any crossings. 

\begin{defn}
A $\textit{virtual knotoid}$ is an equivalence class of virtual knotoid diagrams in $S^2$ up to the virtual equivalence. 
\end{defn}

In \cite{gugumcu2017new}, the first author and Kauffman showed that virtual knotoids in $S^2$ can be considered as knotoids in surfaces as in the case of virtual knots. It is done by extending the stable equivalence relation, defined for knots in surfaces.

\begin{defn}
	Let $K_1$ and $K_2$ be two knotoid diagrams in surfaces $\Sigma_{g_1}$ and $\Sigma_{g_2}$. Then $K_1$ and $K_2$ are called \textit{stably equivalent} if they can be obtained from each other by finitely many Reidemeister moves in the surfaces, homeomorphisms of surfaces and handle stabilization or destabilization in the complement of the diagrams. The resulting equivalence classes are called \textit{stable equivalence classes} of knotoid diagrams in surfaces.
\end{defn}

In \cite{gugumcu2017new}, the authors extended the notion of abstract knot diagrams to abstract knotoid diagrams to provide a one-to-one correspondence between virtual knotoids and stable equivalence classes of knotoid diagrams in surfaces. Abstract knotoid diagrams are defined as follows.

\begin{figure}[h!]
	\centering
	\includegraphics[width=75mm]{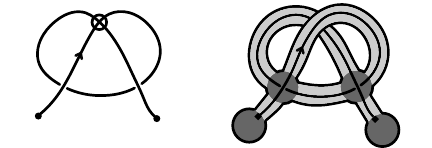}
	\caption{An abstract knotoid diagram obtained from a given virtual knotoid diagram.}
	\label{absknotoid}
\end{figure}

\begin{defn}
    
An \textit{abstract knotoid diagram} is a closed, connected, and orientable ribbon surface with boundary that is associated with a knotoid diagram $K$ as follows.  To each classical crossing and each endpoint of $K$, we attach 2-disks and then connect these disks by ribbons that contain semi-arcs of $K$. Whenever these ribbons meet at a virtual crossing we want them to pass over one another. See Figure \ref{absknotoid} for an illustration.
\end{defn}

Conversely, we can obtain a unique virtual knotoid diagram from an abstract knotoid diagram as follows. Firstly, we embed an abstract knotoid diagram in $S^3$ in such a way that the disks surrounding the classical crossings and the endpoints are mapped on a sphere $S^2$ in $S^3$. Then we project this embedding to $S^2$ in a way such that the projections of the arcs that lie in ribbon bands either do not intersect or intersect transversally under this projection. We mark the transversal intersections in the projection as virtual crossings. The resulting object is a virtual knotoid diagram in $S^2$.

\begin{figure}[h!]
	\centering
	\includegraphics[width=100mm]{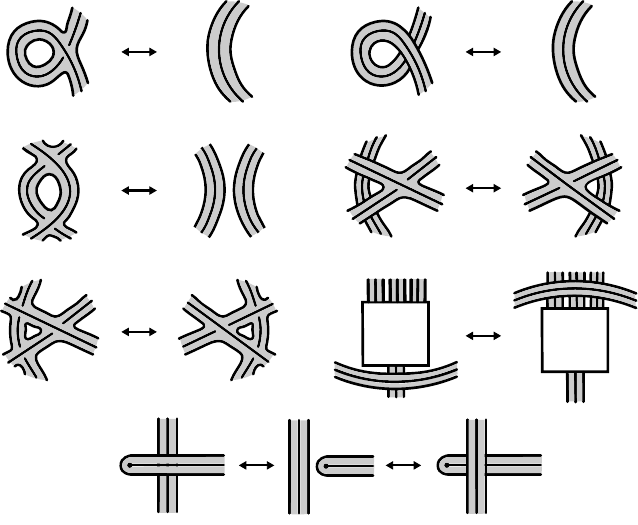}
	\caption{Abstract moves defined on abstract knotoid diagrams.}
	\label{absmoves}
\end{figure}

Abstract knotoid diagrams are considered up to an equivalence relation induced by \textit{generalized abstract moves} depicted in Figure \ref{absmoves}.
An equivalence class of abstract knotoid diagrams is called an \textit{abstract knotoid}.

The two operations given above induce mappings between virtual knotoids in $S^2$ and abstract knotoids, which are inverses of each other. Thus, we have the following theorem.

\begin{thm}
	\cite{gugumcu2017new} There is a one-to-one correspondence between virtual knotoids in $S^2$ and abstract knotoids.
	\label{eq1}
\end{thm}

Let $D$ be an abstract knotoid diagram. By attaching 2-disks to the boundaries of $D$, we obtain a closed, connected and orientable surface containing the knotoid diagram $D$. Conversely, given a knotoid diagram $K$ in the surface $\Sigma_g$, the tubular neighborhood of $K$ is an abstract knotoid diagram. These two operations induce mappings between abstract knotoids and stable equivalence classes of knotoid diagrams. In fact, we have the following theorem.

\begin{thm}
	\cite{gugumcu2017new} There is a one-to-one correspondence between abstract knotoids and stable equivalence classes of knotoid diagrams in surfaces.
	\label{eq2}
\end{thm}

Theorem \ref{eq1} combined with Theorem \ref{eq2} gives the following main theorem.

\begin{thm} \cite{gugumcu2017new}
	There is a one-to-one correspondence between virtual knotoids and stable equivalence classes of knotoid diagrams in surfaces.
	\label{piece1}
\end{thm}

As in the case of virtual knots, classification of virtual knotoids is a challenging problem. This problem has led to the construction of several invariants for virtual knotoids. In \cite{gugumcu2017new}, the first author and Kauffman introduced new invariants of virtual knotoids, such as the parity bracket polynomial, the affine index polynomial, and the arrow polynomial. In \cite{ding2025vassiliev}, Ding et. al. introduced a 0-smoothing polynomial invariant of virtual knotoids and they prove that this invariant is a Vasiliev invariant. In \cite{knotoidTable}, Bartholomew gave a classification table for virtual knotoids of up to seven crossings by utilizing the normalized extended bracket polynomial. 
More recently in \cite{gugumcu2025biquandle}, the authors utilized biquandle virtual brackets to obtain multi-set invariants of virtual knotoids. 

It was conjectured in \cite{gugumcu2017new} that classical knotoid theory properly embeds in virtual knotoid theory, as in the case of classical knots and virtual knots. In the next section, we give a geometric description for virtual knotoids and prove this conjecture.

\section{An interpretation of virtual knotoids in thickened surfaces}
\label{virtualarcs}

In addition to the three-dimensional representation of planar knotoids as line isotopy classes of smooth open curves mentioned in Section \ref{knotoids}, Kodokostas and Lambropoulou introduced \textit{rails arcs} in $\mathbb{R}^3$ as a geometric interpretation of planar knotoids in \cite{kodokostas2019rail}. Here we define rail arcs in thickened surfaces  and consider them up to stable equivalence. Note also that rail arcs in thickened surfaces are special cases of H-curves presented in \cite{chmutov2024thistlethwaite}.

Through the projections of rail arcs, we obtain knotoid diagrams in surfaces, and we establish a correspondence between the equivalence classes of rail arcs and the stable equivalence classes of knotoid diagrams in surfaces.

\subsection{Rail arcs in thickened surfaces}
\label{varcs}
We define a rail arc in a thickened surface as follows.
 
\begin{figure}[h!]
	\centering
	\includegraphics[width=0.37\linewidth]{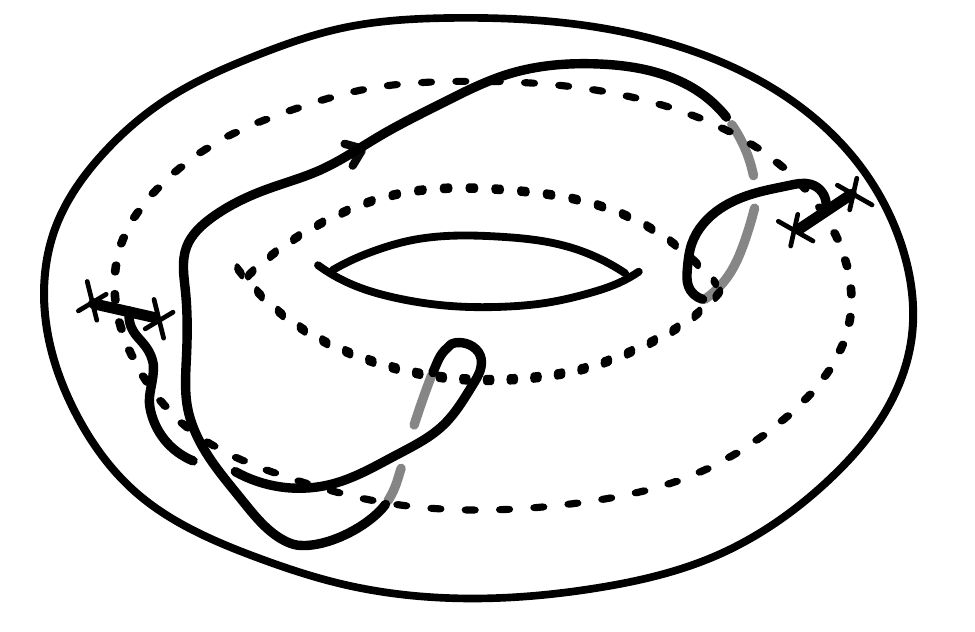}
	\caption{A rail arc in a thickened torus.}
	\label{tthick}
\end{figure}

\begin{defn}
	Let $t$ and $h$ be two distinct points in $\Sigma_g$. A \textit{rail arc} $\boldsymbol{\upsilon}$ is an embedding of the unit interval in the thickened surface $\Sigma_g\times I$ such that $\text{Im}(\boldsymbol{\upsilon})\cap (t\times I)=\boldsymbol{\upsilon}(0)$ and $\text{Im}(\boldsymbol{\upsilon})\cap (h\times I)=\boldsymbol{\upsilon}(1)$. That is, a rail arc is a simple curve in the thickened surface $\Sigma_g\times I$ that starts at a point on $t\times I$ and ends at a point on $h\times I$ (see Figure \ref{tthick} for an example of a rail arc in thickened torus). The points $\boldsymbol{\upsilon}(0)$ and $\boldsymbol{\upsilon}(1)$ that lie on $t\times I$ and $h\times I$, respectively, are the \textit{endpoints} of the arc. The endpoint $\boldsymbol{\upsilon}(0)$ is called the \textit{tail} and the endpoint $\boldsymbol{\upsilon}(1)$ is called the \textit{head} of $\boldsymbol{\upsilon}$. We call the line segments $t\times I$ and $h\times I$ the \textit{rails} of $\boldsymbol{\upsilon}$ following the terminology used in \cite{kodokostas2019rail} in the definition of rail arcs in $\mathbb{R}^3$.
\end{defn}

\begin{defn} Let $\boldsymbol{\upsilon}$ be a rail arc in the thickened surface $\Sigma_g \times I$ with rails $t\times I$ and $h\times I$. $\boldsymbol{\upsilon}$ is called \textit{trivial} if it embeds in $\sigma\times I$ where $\sigma$ is a rail arc in the surface $\Sigma_g\times \{0\}$ with endpoints $t\times\{0\}$ and $h\times \{0\}$. 
\end{defn}

We now introduce an analog of the isotopy given in Definition \ref{li} for thickened surfaces. Then we describe the handle stabilization and destabilization operations for a rail arc in a thickened surface.

\begin{defn}
	Let $\boldsymbol{\upsilon}_i$ be a rail arc in $\Sigma_g\times I$ with its rails $t_i\times I$ and $h_i\times I$, for $i=1,2$. $\boldsymbol{\upsilon}_1$ and $\boldsymbol{\upsilon}_2$ are said to be \textit{arc isotopic} if there is a smooth ambient isotopy of $\Sigma_g\times I$, taking one arc to the other arc in the complement of the rails $t_i\times I$ and $h_i\times I$, taking endpoints to endpoints; $(t_1,x_1)$ to $(t_2,x_2)$ and $(h_1,y_1)$ to $(h_2,y_2)$, and taking rails to rails; $t_1 \times I$ to $t_2 \times I$ and $h_1 \times I$ to $h_2 \times I$.
	\label{def}
\end{defn}

\begin{rem}
    \label{forbid}
    Note that, since the ambient isotopy that takes one arc to the other is performed in the complement of the rails of the rail arc in Definition \ref{def}, a rail arc does not pass through its rails during the arc isotopy.
\end{rem}

In our interpretation, we need to take extra care of the rails $t\times I$ and $h\times I$ of the rail arc $\boldsymbol{\upsilon}\subset \Sigma_g\times I$. In particular, whenever we perform an operation, aside from the ambient isotopy in the ambient thickened surface, we need to make sure that we work in the complement of not only the rail arc itself, but also its rails. In other words, instead of $(\Sigma_g\times I) \;-\; \boldsymbol{\upsilon}$, the operation is performed in $(\Sigma_g\times I)\;-\;(\boldsymbol{\upsilon}\cup (\{t,h\}\times I))$.

\begin{defn}
	Consider a rail arc $\boldsymbol{\upsilon}$ in a thickened surface $\Sigma_g \times I$ with its rails $t\times I$ and $h\times I$. Let $\boldsymbol{|\upsilon|}$ denote the union of $\boldsymbol{\upsilon}$ and its rails $t\times I$ and $h\times I$. Let $B_1, B_2$ be two disjoint open disks in $\Sigma_g$ such that  $B_1\times I$ and $B_2\times I$ does not intersect $\boldsymbol{|\upsilon|}$. The \textit{stabilization} operation consists of removing the two open cylinders $B_1\times I$ and $B_2\times I$ from $\Sigma_g\times I$, and attaching a thickened handle $S^1\times I \times I$ by gluing annuli $S^1\times \{0\}\times I$ and $S^1\times \{1\}\times I$ to the boundaries of removed open cylinders in $\Sigma_g\times I$. Let $C$ be a simple closed curve in $\Sigma_g$ such that $C\times I$ is disjoint from $\boldsymbol{|\upsilon|}$ and it does not bound a ball in $(\Sigma_g\times I)-\boldsymbol{|\upsilon|}$. The \textit{destabilization} operation consists of cutting $\Sigma_g\times I$ along the vertical annulus $C\times I$, and capping the two resulting annuli with thickened disks $D^2\times I$.
	\label{destabknotoid}
\end{defn}

\begin{rem}
	We consider destabilization annuli up to isotopy in $\Sigma_g \times I$, following the same notion as in \cite{kuperberg2003virtual}. We let any isotopic image of the annulus $C\times I$ mentioned in Definition \ref{destabknotoid} to be a destabilization annulus.
\end{rem}

\begin{defn}
	The arc isotopies of $\;\Sigma_g\times I$, homeomorphisms of $\Sigma_g$ and stabilizations and destabilizations of  $\Sigma_g\times I$ together induce an equivalence relation on the set of rail arcs in thickened surfaces. We call this equivalence relation \textit{stable equivalence}. Two rail arcs that are equivalent up to rail arc equivalence are called \textit{stably} equivalent.
\end{defn}

From now on, we consider equivalence classes of rail arcs with respect to the stable equivalence. We use the notation $(\boldsymbol{\upsilon},\Sigma_g\times I)$ for a rail arc $\boldsymbol{\upsilon}$ in the thickened surface $\Sigma_g\times I$. Let $t\times I$ and $h\times I$ be the rails of $\boldsymbol{\upsilon}$. The projection of $\boldsymbol{\upsilon}$ and its rails on the surface $\Sigma_g\times \{0\}$, which we identify with $\Sigma_g$, by the map 
    \begin{align}
    \label{proj}
    	p:\Sigma_g\times I&\rightarrow \Sigma_g\\
    	(x,t)&\mapsto x\notag
    \end{align}
    equiped with under or over information at the crossings gives a diagram of $\boldsymbol{\upsilon}$. Observe that this type of a projection does not cause any virtual crossings. Hence, the resulting diagram is a classical knotoid diagram $K$ in the surface $\Sigma_g$ with tail $t\times\{0\}$ and head $h\times\{0\}$, which we denote by $(K,\Sigma_g)$.

\begin{figure}[h!]
	\centering
	\includegraphics[width=60mm]{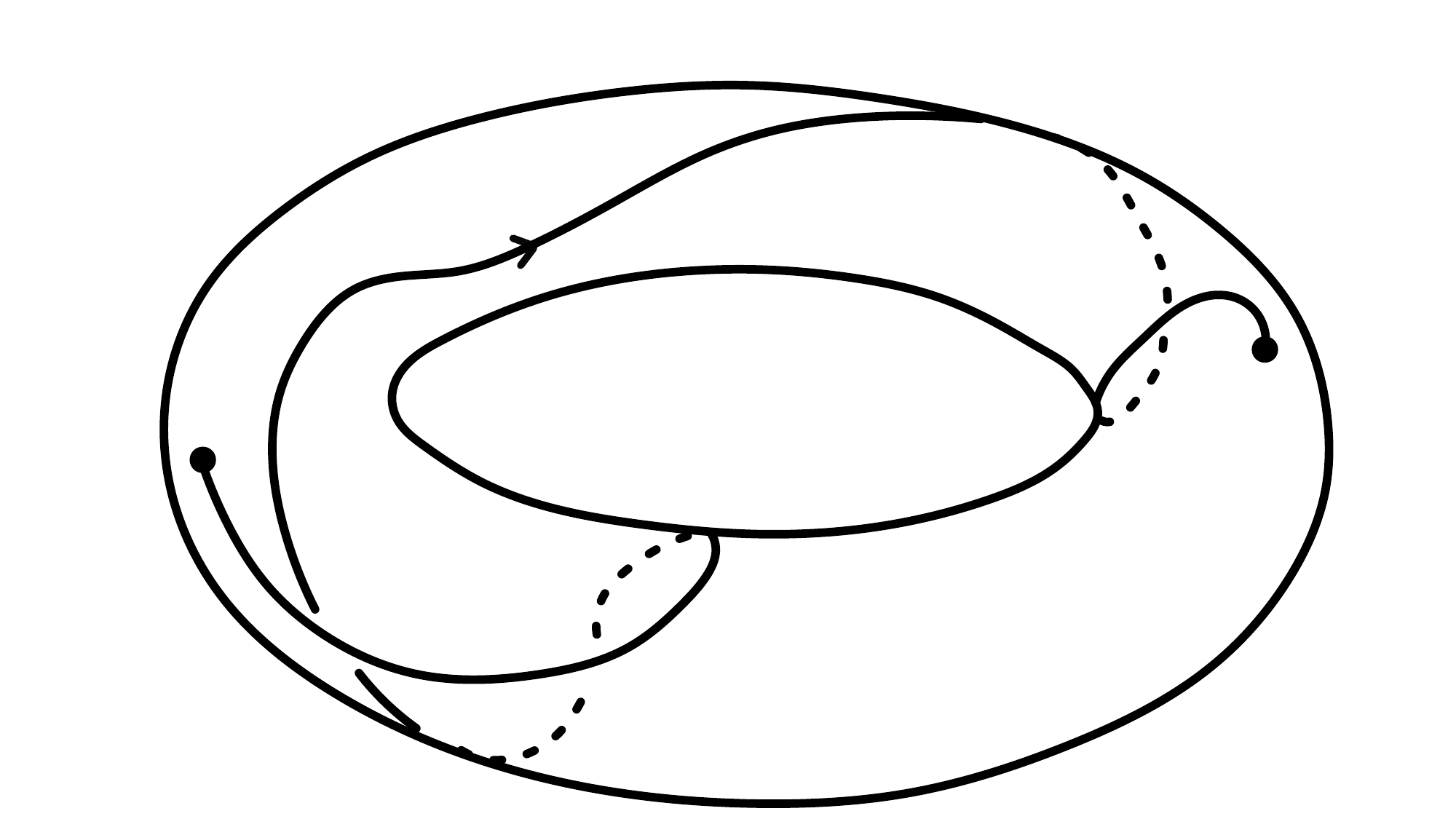}
	\caption{The diagram of the rail arc given in Figure \ref{tthick}.}
	\label{torthick}
\end{figure}

\begin{prop}
    \label{stabs}
    Let $(\boldsymbol{\upsilon},\Sigma_g\times I)$ be a rail arc and $(K,\Sigma_g)$ be its diagram. The handle stabilizations and destabilizations of $(\boldsymbol{\upsilon},\Sigma_g\times I)$ correspond to handle stabilizations and destabilizations of $(K,\Sigma_g)$, and vice versa.
\end{prop}

\begin{proof}
    Consider the open disk $B_{i=1,2}$, and the simple closed curve $C$ in the definition of stabilization in Definition \ref{destabknotoid}. Since $B_i\times I$ does not intersect $\boldsymbol{|\upsilon|}$, $B_i$ does not intersect the diagram $K$ in $\Sigma_g$. Hence, the operation that consists of removing $B_i\times I$ from $\Sigma_g\times I$ and attaching a thickened handle $S^1\times I \times I$ to $\Sigma_g\times I$ removes the open disk $B_i$ from $\Sigma_g$ and attaches a handle $S^1\times I$ to $\Sigma_g$ in the complement of $K$, which corresponds to a stabilization of $(K,\Sigma_g)$. Similarly, since $C\times I$ does not intersect $\boldsymbol{|\upsilon|}$, $C$ does not intersect $K$. Morever, since $C\times I$ does not bound a ball in the complement of $\boldsymbol{|\upsilon|}$, $C$ does not bound a disk in the complement of $K$ in $\Sigma_g$. Hence, destabilization operation on $(\boldsymbol{\upsilon},\Sigma_g\times I)$ induces a destabilization on $(K,\Sigma_g)$. As the diagram $(K,\Sigma_g)$ of $(\boldsymbol{\upsilon},\Sigma_g\times I)$ is obtained through the projection \eqref{proj}, the other direction of the proof follows using similar arguments, thus omitted here.
\end{proof}

\begin{thm}
	Two rail arcs $(\boldsymbol{\upsilon},\Sigma_{g}\times I)$ and $(\boldsymbol{\upsilon}',\Sigma_{g'}\times I)$ are stably equivalent if and only if their diagrams $(K,\Sigma_{g})$ and $(K',\Sigma_{g'})$ are stably equivalent.
	\label{equi}
\end{thm}
To prove Theorem \ref{equi}, we utilize piecewise-linear representations of rail arcs since rail arcs are in the smooth category. A piecewise-linear representation of a rail arc is the union of finitely many edges $[p_1p_2],\; [p_2p_3],\;...,\;[p_{n-1}p_n]$, such that $p_1$ and $p_n$ correspond to the tail and head of the rail arc, respectively. Each pair of edges in the representation of a rail arc can only intersect at a point $p_i$, where $i\in\{2,3,...,n-1\}$.

\begin{defn}
	Let $[p_ip_{i+1}]$ be an edge of a piecewise-linear representation of a rail arc $\boldsymbol{\upsilon}$ in $\Sigma_g\times I$, and let $p_0$ be a point in the complement of $\boldsymbol{|\upsilon|}$, $(\Sigma_g\times I)- \boldsymbol{|\upsilon|}$, such that $[p_ip_0],[p_0p_{i+1}]\subset \Sigma_g\times I$ and the triangle $\overset{\triangle}{p_ip_0p_{i+1}}$ bounds a disk in $(\Sigma_g\times I)- \boldsymbol{|\upsilon|}$. In other words, the triangular region $\overset{\triangle}{p_ip_0p_{i+1}}$ intersects neither the rail arc $\boldsymbol{\upsilon}$ nor its rails. We define the \textit{triangle move} in the thickened surface $\Sigma_g\times I$ as the transformation of $[p_ip_{i+1}]$ into two edges $[p_ip_0]$ and $[p_0p_{i+1}]$. The \textit{inverse triangle move} is the transformation of two consecutive edges into one edge that satisfies the condition that the triangular region formed by these three edges bounds a disk in $\Sigma_g\times I-\boldsymbol{|\upsilon|}$. These two moves are denoted by $\Delta^{+1}$ and $\Delta^{-1}$, respectively.
\end{defn}

The equivalence of rail arcs can be expressed as the equivalence of piecewise-linear representations of rail arcs. Two representations are equivalent if one can be deformed into the other by a finite sequence of $\Delta^{\pm 1}$ moves.
Equipped with these, we now give the proof of Theorem \ref{equi}.

\begin{figure}[h!]
	\centering
	\includegraphics[width=80 mm]{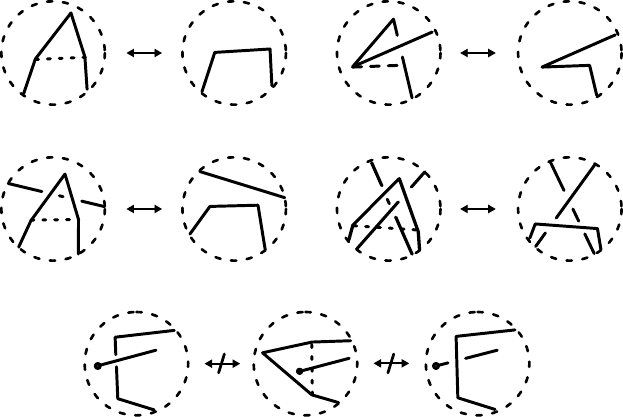}
	\caption{Shadows of triangle moves in surface $\Sigma$.}
	\label{shadows}
\end{figure}

\begin{proof} (of Theorem \ref{equi}) Assume that $(\boldsymbol{\upsilon},\Sigma_g\times I)$ and $(\boldsymbol{\upsilon}',\Sigma_{g'}\times I)$ are stably equivalent.
	First, note that we can consider the steps of obtaining $(\boldsymbol{\upsilon}',\Sigma_{g'}\times I)$ from $(\boldsymbol{\upsilon},\Sigma_g\times I)$ in two parts; handle stabilizations or destabilizations and arc isotopies.
    
    Let $(\boldsymbol{\upsilon}'',\Sigma_{g''}\times I)$ be a rail arc obtained by applying a handle stabilization or destabilization on $(\boldsymbol{\upsilon},\Sigma_g\times I)$. By Proposition \ref{stabs}, this handle stabilization or destabilization performed on $(\boldsymbol{\upsilon},\Sigma_g\times I)$ is captured by a handle stabilization or destabilization of $(K,\Sigma_g)$, which transforms $(K,\Sigma_g)$ into the diagram $(K'',\Sigma_{g''})$ of $(\boldsymbol{\upsilon}'',\Sigma_{g''}\times I)$.

    To prove that arc isotopies of $(\boldsymbol{\upsilon},\Sigma_g\times I)$ are captured by isotopies of $\Sigma_g$ and Reidemeister moves on $(K,\Sigma_g)$, we first think of these arc isotopies as series of triangle moves $\Delta^{\pm 1}$ on the piecewise-linear representation of $(\boldsymbol{\upsilon},\Sigma_g\times I)$. Then we show that all possible shadows of the triangle moves are generated by isotopies of $\Sigma_g$ and Reidemeister moves on $(K,\Sigma_g)$. We use induction on the number of strands, denoted by $n$, in the shadow of a triangle move. Figure \ref{shadows} shows that the isotopies of $\Sigma_g$ and Reidemeister moves are, respectively, enough for the cases $n=0,1,2$ where the lines intersect for $n=2$. For non-intersecting lines in the case $n=2$, and for $n\geq 3$, we apply subdivision on the edges of the triangle and express the triangle move as a finite sequence of triangle moves of types given in Figure \ref{shadows}.

    Moreover, remember that in a triangle move, the triangle does not intersect a rail. In other words, during a triangle move, no parts of the piecewise-linear representation of a rail arc pass through a rail. This means that the shadow of the triangle of a triangle move does not contain an endpoint of the diagram. Hence, the realizations of triangle moves on the surface do not violate forbidden moves for knotoid diagrams as shown in Figure \ref{shadows} as well. This shows that an arc isotopy that transforms $(\boldsymbol{\upsilon},\Sigma_g\times I)$ into $(\boldsymbol{\upsilon}'',\Sigma_{g''}\times I)$ is captured by isotopies of $\Sigma_g$ and Reidemeister moves, which transforms $(K,\Sigma_g)$ into $(K'',\Sigma_{g''})$ without violating forbidden moves.
    
    We conclude that handle stabilizations or destabilizations and arc isotopies of thickened surfaces that transform $(\boldsymbol{\upsilon},\Sigma_g\times I)$ into $(\boldsymbol{\upsilon}',\Sigma_{g'}\times I)$ are captured by handle stabilizations or destabilizations and isotopies of surfaces together with Reidemeister moves that transform $(K,\Sigma_g)$ into $(K',\Sigma_{g'})$. Hence, stable equivalence of $(\boldsymbol{\upsilon},\Sigma_g\times I)$ and $(\boldsymbol{\upsilon}',\Sigma_{g'}\times I)$ implies stable equivalence of their diagrams $(K,\Sigma_g)$ and $(K',\Sigma_{g'})$. The other direction follows from similar arguments and Proposition \ref{stabs} as well and is omitted.
\end{proof}

\subsection{Unifying the interpretations}
\label{main}

As explained in Section~\ref{vknotoids}, virtual knotoid theory has an interpretation as stable equivalence classes of knotoid diagrams in surfaces. An interpretation of the latter theory is given in Section~\ref{varcs} as rail arcs in thickened surfaces. Here we combine these two interpretations to give a geometric interpretation of virtual knotoids.

\begin{thm}
	The rail arc theory in thickened surfaces is equivalent to the virtual knotoid theory given diagrammatically in $S^2$.
	\label{vct}
\end{thm}

\begin{proof}
	Theorem \ref{piece1} states that virtual knotoids in $S^2$ are in one-to-one correspondence with the stable equivalence classes of knotoid diagrams in surfaces. Theorem \ref{equi} states that stable equivalence classes of rail arcs in thickened surfaces are in one-to-one correspondence with stable equivalence classes of knotoid diagrams in the underlying surfaces. By combining these two theorems, we obtain a correspondence between virtual knotoids in $S^2$ and stable equivalence classes of rail arcs in thickened surfaces.
\end{proof}

Note that Theorem \ref{vct} allows us to study virtual knotoids in $S^2$ through stable equivalence classes of rail arcs in thickened surfaces.

\begin{rem}
	By fixing the thickened surface as $S^2\times I$, and restricting the equivalence relation to be up to only arc isotopies of $S^2\times I$ in our interpretation, we obtain a geometric interpretation of classical knotoids in $S^2$. A classical knotoid in $S^2$ may also be regarded as a virtual knotoid when considered up to virtual equivalence. In this setting, a classical knotoid is a virtual knotoid with at least one irreducible representation in $S^2 \times I$.
\end{rem}

\begin{defn}
	Let $(\boldsymbol{\upsilon},\Sigma_g\times I)$ be a rail arc. If the complement of $\boldsymbol{|\upsilon|}$ in $\Sigma_g\times I$ does not have empty thickened handles, then $(\boldsymbol{\upsilon},\Sigma_g\times I)$ is called an \textit{irreducible rail arc}.
\end{defn}

Given a rail arc in a thickened surface, we can remove the empty handles by the destabilization operation and represent this rail arc by an irreducible rail arc. It is natural to ask if this representation is unique up to arc isotopies. 

We use Kuperberg's arguments in \cite{kuperberg2003virtual} in our setting which considers rail arcs to prove that every rail arc admits a unique representation by an irreducible rail arc. Before we state the theorem, we point out some differences that arise when we consider rail arcs in thickened surfaces.
    
Links in thickened surfaces are multi-component objects. Kuperberg extended the destabilization operation by allowing a split component of a link in a thickened surface $\Sigma_g\times I$, which is separated from the rest of the link by a sphere or a disk, to be removed from $\Sigma_g\times I$ and placed in a separate $S^2\times I$. For this reason, the surfaces allowed in his definition of destabilization were vertical annuli, spheres, and proper disks, which he called \textit{admissible}. He also required the surface to be \textit{essential} in the sense that the surface itself, or together with a boundary component of the thickened surface, does not bound a ball in the complement of the link.

\begin{prop}
\label{spc}
    Let $(\boldsymbol{\upsilon},\Sigma_g\times I)$ be a rail arc. Then, spheres and properly embedded disks in $\Sigma_g\times I-\boldsymbol{|\upsilon|}$ are inessential.
\end{prop}
\begin{proof}
    By definition, a rail arc with its rails is of a single component. Thus, a sphere is essential in $\Sigma_g\times I$ if and only if it bounds an open ball $B$ that contains the rail arc together with its rails, $\boldsymbol{|\upsilon|}$. Similarly, a properly embedded disk is essential in $\Sigma_g\times I$ if and only if the disk and a part of the boundary of $\Sigma_g\times I$ bound an open ball $B$ that contains $\boldsymbol{|\upsilon|}$. However, rails are line segments between the inner boundary component $\Sigma_g\times \{0\}$ and the outer boundary component $\Sigma_g\times \{1\}$. In other words, if we denote the rails of $\boldsymbol{\upsilon}$ by $h\times I$ and $t\times I$, then the endpoints $\{h,t\}\times \{i\}$ are in $\Sigma_g\times \{i\}$ for $i=0,1$. This means that these four endpoints of the rails are not contained in any open ball in $\Sigma_g\times I$. Hence, such open ball $B$ does not exist in $\Sigma_g\times I$.
\end{proof}
    
As a result of Proposition \ref{spc}, the only possible admissible and essential surfaces are properly embedded vertical annuli in the complement of the rail arc and its rails in our setting.

A simple closed curve $C$ in $\Sigma_g$ is \textit{essential} if it does not bound a disk in $\Sigma_g$. In the case of links in thickened surfaces, an inessential curve $C$ can give rise to an essential annulus $A=C\times I$. This happens when $D\times I$, where $D$ is the disk bounded by $C$, contains a split component of the link. In the case of rail arcs in thickened surfaces, $C$ is essential if it bounds a disk  in the complement of $\{h,t\}\times \{0\}$ in $\Sigma_g\times \{0\}$. The following proposition states that an inessential simple closed curve cannot give rise to an essential annulus in the case of rail arcs in thickened surfaces.

\begin{prop}
    \label{spc1}
    Let $(\boldsymbol{\upsilon},\Sigma_g\times I)$ be a rail arc with rails $\{h,t\}\times I$. Let $A=C \times I$ be a vertical annulus in the complement of $\boldsymbol{|\upsilon|}$ in $\Sigma_g\times I$ where $C$ is a closed curve in $\Sigma_g\times \{0\}$. $C$ is essential if and only if $A$ is essential.
\end{prop}
\begin{proof}
    Let $C$ be inessential. Then, $C$ bounds a disk $D$ which does not contain the endpoints $\{h,t\}\times \{0\}$ in $\Sigma_g\times \{0\}$. Since $\boldsymbol{|\upsilon|}$ is of a single component, $(D\times I)\cap \boldsymbol{|\upsilon|}=\emptyset$. Hence, $A$ bounds a ball in the complement of $\boldsymbol{|\upsilon|}$, which means $A$ is inessential. The other direction follows by reversing the arguments.
\end{proof}

By Propositions \ref{spc} and \ref{spc1}, the destabilization operation on rail arcs is performed along vertical annuli $A=C\times I$, where $C$ is essential in $\Sigma_g$. Hence, we refer such annuli as \textit{destabilization annuli}. We use the term destabilization to refer both to the operation and the resulting thickened surface from applying the operation along a destabilization annulus.

\begin{figure}[h!]
    \centering
    \includegraphics[width=0.5\linewidth]{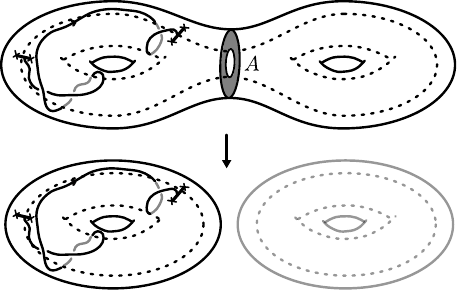}
    \caption{Destabilization along a separating annulus $A$.}
    \label{sdesan}
\end{figure}

Let $(\boldsymbol{\upsilon},\Sigma_g\times I)$ be a rail arc. A destabilization annulus $A$ in $\,\Sigma_g\times I\,$ is called \textit{separating} if $\Sigma_g\times I- A$ is disconnected and \textit{non-separating} otherwise. It is clear that destabilization along a non-separating annulus $A$ removes an empty thickened handle from $\Sigma_g\times I$, hence reduces the genus by one. If $A$ is separating, then destabilization along $A$ splits $\Sigma_g\times I$ into two thickened surfaces $\Sigma_{g_1}\times I$ and $\Sigma_{g_2}\times I$, where $g=g_1+g_2$. Since $\boldsymbol{|\upsilon|}$ is of a single component, it is contained in exactly one of these thickened surfaces, say $\Sigma_{g_1}\times I$. In this case, we discard $\Sigma_{g_2}\times I$ whose intersection with $\boldsymbol{|\upsilon|}$ is empty. Hence, destabilization along $A$ decreases the genus of the thickened surface in which $\boldsymbol{\upsilon}$ lies exactly by $g_2$. See Figure \ref{sdesan} for an illustration. In this illustration, we see the destabilization of a rail arc in thickened double torus along a separating annulus $A$. As a result, we obtain two thickened tori. The one on the left-hand side contains the rail arc. The one on the right-hand side does not intersect the rail arc and its rails. Hence, we discard the thickened torus on the right-hand side.

\begin{thm}
	Every rail arc is uniquely represented by an irreducible rail arc.
	\label{thm1}
\end{thm}

\begin{proof}
    Assume to the contrary that there exist rail arcs that admit multiple irreducible representations. Let $\Sigma_g \times I$ be the least genus thickened surface among the thickened surfaces in which such rail arcs can lie in and $\mathcal{M} = \{(\boldsymbol{\upsilon},\Sigma_g\times I)\}$ be the collection of rail arcs in $\Sigma_{g} \times I$ that admit multiple irreducible representations. Then, for a rail arc in $\mathcal{M}$,  there exist different destabilization annuli so that we may end up with multiple destabilizations that are not arc isotopic. Any destabilization of $(\boldsymbol{\upsilon},\Sigma_g\times I)$ has genus less than $g$, and by assumption, each of these destabilizations admits a unique irreducible representation. We show that the unique irreducible representations of destabilizations of $(\boldsymbol{\upsilon},\Sigma_g\times I)$ are arc isotopic by analyzing the relative positions of their destabilization annuli.
    	
    Let $A_1$ and $A_2$ be two disjoint annuli along which we can destabilize $(\boldsymbol{\upsilon},\Sigma_g\times I) \in \mathcal{M}$. Let $(\boldsymbol{\upsilon},\Sigma_{g_1}\times I)_{A_1}$ and $(\boldsymbol{\upsilon},\Sigma_{g_2}\times I)_{A_2}$ denote destabilizations of $(\boldsymbol{\upsilon},\Sigma_g\times I)$ along $A_1$ and $A_2$. If $A_1$ and $A_2$ are parallel in $\Sigma_g\times I$, then it is clear that the destabilizations along these annuli are arc isotopic. If they are not parallel, we destabilize $(\boldsymbol{\upsilon},\Sigma_{g_1}\times I)_{A_1}$ and $(\boldsymbol{\upsilon},\Sigma_{g_2}\times I)_{A_2}$ along $A_2$ and $A_1$, respectively, to obtain a common destabilization $(\boldsymbol{\upsilon},\Sigma_{g_c}\times I)_C$. Notice that $(\boldsymbol{\upsilon},\Sigma_{g_c}\times I)_C$ has a unique irreducible representation which is also the unique irreducible representation of both $(\boldsymbol{\upsilon},\Sigma_{g_1}\times I)_{A_1}$ and $(\boldsymbol{\upsilon},\Sigma_{g_2}\times I)_{A_2}$ since each of them has genus less than $g$. We illustrate this in Figure \ref{fig:des}. This sets a contradiction with the assumption that $(\boldsymbol{\upsilon},\Sigma_g\times I)$ has multiple inequivalent irreducible representations.
    \begin{figure}[h!]
        \centering
        \includegraphics[width=0.55\linewidth]{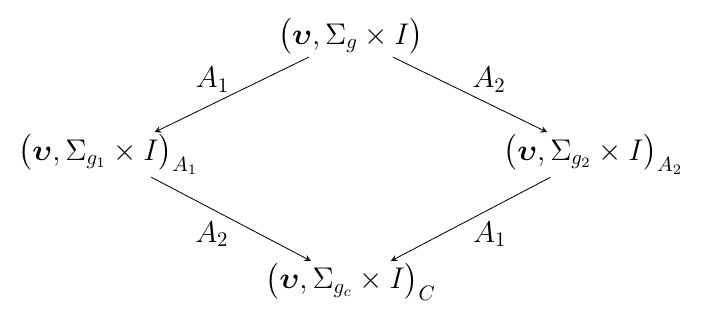}
        \vspace{-0.2 cm}
        \caption{A diagram of destabilizations $(\boldsymbol{\upsilon},\Sigma_g\times I)$ along annuli $A_1$ and $A_2$. An arrow with an indicator $A_i$ means that the rail arc at the head of the arrow is obtained by destabilization of the rail arc at the tail of the arrow along $A_i$.}
        \label{fig:des}
    \end{figure}
    	
    Now we consider non-disjoint pairs of destabilization annuli of $(\boldsymbol{\upsilon},\Sigma_g\times I)$ such that the destabilizations of $(\boldsymbol{\upsilon},\Sigma_g\times I)$ along these annuli have inequivalent unique irreducible representations. Among such pairs of annuli, consider $A_1$ and $A_2$ in general position so that they intersect transversely at finitely many disjoint curves, and the number of intersection curves is minimal. Let the number of curves in $A_{1} \cap A_{2}$ be $k$. Note that this assumption implies that any two destabilizations of $(\boldsymbol{\upsilon},\Sigma_g\times I)$ along a pair of annuli that intersect at less than $k$ curves have equivalent unique irreducible representations.
    
    Let $C$ be a curve in $A_1 \cap A_2$. Then, there are four possible cases. If $C$ is open, either it is \textit{vertical} so that its endpoints are in different components of $\partial A_{i=1,2}$, see Figure \ref{intarc}, or it bounds a disk in $A_i$ with a part of $\partial A_i$. If $C$ is closed, either it is  \textit{horizontal} so that it is parallel to the boundaries of $A_1$ and $A_2$, see Figure \ref{intc}, or it bounds a disk in $A_i$.
    
    \begin{figure}[h!]
        \centering    
        \begin{subfigure}{0.9\textwidth}
            \centering
            \includegraphics[width=55mm]{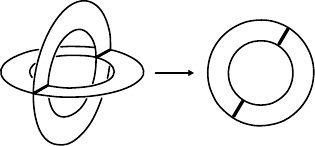}
            \vspace{0.3 cm}
            \caption{Vertical intersection open curves.}
            \label{intarc}
        \end{subfigure}\hfill
        \begin{subfigure}
            {0.9\textwidth}
            \vspace{0.2cm}
            \centering
            \includegraphics[width=46mm]{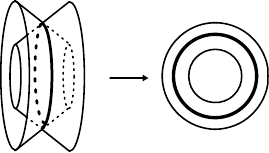}
            \caption{A horizontal intersection closed curve.}
            \label{intc}
        \end{subfigure}
        \caption{Examples of intersection curves.}
        \label{iac}
    \end{figure}

    We first consider the case where $C$ is a non-horizontal closed curve or a non-vertical open curve in $A_i$. Then, $C$ is either the boundary or part of the boundary of a disk, say $D$, in $A_i$. We refer $D$ as the disk associated to $C$. $C$ is called \textit{innermost} if $D$ does not contain any other intersection curves. Since $A_1\cap A_2$ consists of finitely many, disjoint intersection curves, if $A_1\cap A_2$ has such a curve $C$, then there exist an innermost intersection curve $C'$ in $A_i$. Let $D'$ be the associated disk of $C'$.   
    We compress $A_{3-i}$ along $D'$ (see Figure \ref{comp} for an illustration of such compression). This operation splits $A_{3-i}$ into two disjoint parts: a destabilization annulus $A_{3-i}^{'}$ and a part that bounds a ball in the complement of $\boldsymbol{|\upsilon|}$. The intersections $A_i\cap A_{3-i}^{'}$ and $A_{3-i}\cap A_{3-i}^{'}$ both consist of fewer than $k$ curves. Thus, the destabilizations along $A_i$, $A_{3-i}$, and $A_{3-i}^{'}$ have equivalent unique irreducible representations. This sets a contradiction.

	\begin{figure}[h!]
		\centering
		\includegraphics[width=0.35\linewidth]{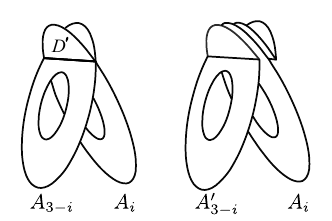}
		\vspace{0.2 cm}
		\caption{An illustration of the compression.}
		\label{comp}
	\end{figure}
	
	Now suppose that $A_1\cap A_2$ consists of $k$ vertical open curves. Consider $A_1\cup A_2$.
	Let $N$ be a regular neighborhood of $A_1\cup A_2$ that does not intersect $\boldsymbol{|\upsilon|}$. Then $\partial N=U_1\cup U_2 \cup ... \cup U_n$, where $U_j$ is a vertical annulus in $\Sigma_g\times I$ for $j=1,...,n$ (see Figure \ref{lastcase} for an example case where $n=4$). Clearly, $A_i\cap U_j=\emptyset$ for all $j =1, ...n$, and $\Sigma_g\times I \;-\; \partial N$ consists of a number of disjoint components. Notice that $\boldsymbol{|\upsilon|}$ lies in exactly one of these components. Now we examine the possible cases for $U_{j=1,...,n}$.
	
	\begin{itemize}
		\item Suppose that none of $U_{j=1,...,n}$ is a destabilization annulus. That is, each of them bounds a ball in $\Sigma_g\times I-\boldsymbol{|\upsilon|}$. Since the regular neighborhood $N$ is surrounded by these annuli $U_{j=1,...,n}$, then one of them, say $U_i$, bounds a ball $B$ that completely contains $N$ in $\Sigma_g\times I-\boldsymbol{|\upsilon|}$. It means that $A_1$ and $A_2$ are contained in $B$, as well. Since $A_1$ is a vertical annulus, the boundaries of $A_1$ are in $\partial B$, and hence bound disks $D_0$ and $D_1$ in $\partial B$. Then the union $D_0\cup A_1\cup D_1$ forms a sphere that bounds a ball in the complement of $\boldsymbol{|\upsilon|}$. Hence, $A_1$ is inessential, which is a contradiction.
        
		\item Suppose now that for some $k\in\{1,...n\}$, $U_k$ is a destabilization annulus. Since $U_k$ does not intersect $A_1$ or $A_2$, destabilization along $U_k$ induces an equivalent unique irreducible representation with destabilization along $A_1$ and $A_2$, by assumption. Thus, the destabilizations of $(\boldsymbol{\upsilon},\Sigma_g\times I)$ along $A_1$ and $A_2$ have equivalent unique irreducible representations, which is a contradiction.
	\end{itemize}
	\begin{figure}[h!]
		\centering
		\includegraphics[width=0.24\linewidth]{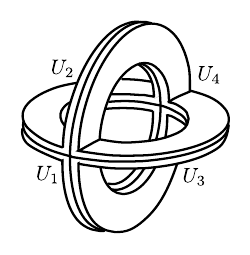}
		\caption{The boundary of a regular neighborhood of $A_1\cup A_2$, where $A_{i=1,2}$ is a destabilization annulus. }
		\label{lastcase}
	\end{figure}
	
	Next, suppose that $A_1\cap A_2$ consists of closed horizontal curves $C_1,...,C_k$, that is, $C_{j=1,...,k}$ does not bound a disk in $A_{i=1,2}$. Let the indexing be such that $C_1$ is the first curve and $C_k$ is the last curve we encounter as we move from the inner boundary ($\partial_0A_i$) to the outer boundary ($\partial_1A_i$) of $A_i$, for each $i=1,2$. Then, $\partial_0A_1$ and $C_1$ bound an annulus $A_1'\subset A_1$ such that $A_1'\cap A_2=C_1$. $C_1$ and $\partial_1A_2$ also bound an annulus $A_2'\subset A_2$ which may have an intersection with $A_1$ other than $C_1$. The union $A_1'\cup A_2'$ is a vertical annulus in $\Sigma_g\times I$ with boundary curves $\partial_0A_1$ and $\partial_1A_2$. Since $A_1$ and $A_2$ are essential, it follows from Proposition \ref{spc1} that their boundary curves $\partial_iA_1$ and $\partial_iA_2$ are essential in $\Sigma_g\times \{i\}$ for $i=0,1$. Hence, $A_1'\cup A_2'$ has essential boundary curves, which means $A_1'\cup A_2'$ is essential and thus it is a destabilization annulus. See Figure \ref{newborn} for an illustration. The number of intersection curves in $(A_1'\cup A_2')\cap A_1$ and $(A_1'\cup A_2')\cap A_2$ is less than $k$. Thus, destabilizations along $A_1'\cup A_2'$, $A_1$, and $A_2$ have equivalent unique irreducible representations, which again is a contradiction.
	\begin{figure}[h!]
		\centering
		\includegraphics[width=0.53\linewidth]{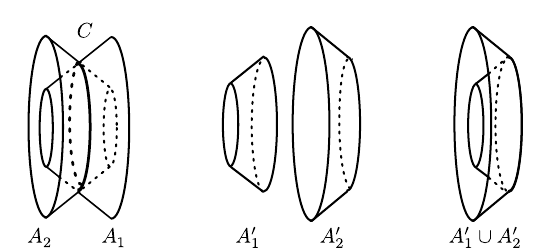}
		\vspace{0.3 cm}
		\caption{Formation of a destabilization annulus $A_1'\cup A_2'$ through a horizontal intersection curve $C$ in the intersection of two destabilization annuli $A_1$ and $A_2$.}
		\label{newborn}
	\end{figure}
	
	The existence of annuli along which the resulting destabilizations have non arc isotopic unique irreducible representations leads to contradictions in all possible cases. Thus, we conclude that a rail arc $(\boldsymbol{\upsilon},\Sigma_g\times I)$ with multiple inequivalent irreducible representations does not exist.
\end{proof}

\begin{thm}
	If two rail arcs in $S^2\times I$ are stably equivalent (up to arc isotopy and handle stabilizations or destabilizations) to each other, then they are equivalent to each other in $S^2\times I$ with respect to arc isotopy.
	\label{clinv}
\end{thm}
\begin{proof} Let $\boldsymbol{\upsilon}_1$ and $\boldsymbol{\upsilon}_2$ be two stably equivalent rail arcs in $S^2\times I$. Then, $\boldsymbol{\upsilon}_1$ and $\boldsymbol{\upsilon}_2$ are equivalent up to arc isotopy in $S^2\times I$ by Theorem \ref{thm1}, since they are already irreducible in $S^2\times I$.
\end{proof} 

By the correspondence given in Theorem \ref{vct} and Theorem \ref{clinv} we have the following corollary which was in fact given as a conjecture in \cite{gugumcu2017new}.

\begin{cor}
	 If two classical knotoid diagrams in $S^2$ are equivalent to each other up to generalized Reidemeister moves, then they are equivalent to each other up to classical Reidemeister moves.
\end{cor}

\end{document}